\input amstex
\documentstyle{amsppt}
\magnification=\magstep1

\pageheight{9.0truein}
\pagewidth{6.5truein}

\TagsOnRight

\hyphenation{uni-ser-ial}

\input xy
\xyoption{matrix}\xyoption{arrow}

\long\def\ignore#1{#1}

\def\la{{\Lambda}}
\def\lamod{\Lambda\text{-}\roman{mod}}
\def \len{\operatorname{length}} 
\def\detour{\mathrel{\wr\wr}}

\def\SS{{\Bbb S}}
\def\Hom{\operatorname{Hom}}

\def\GL{\operatorname{GL}}
 
\def\G{{\Cal G}}

\def\rep{\operatorname{Rep}}
\def\modla#1{{\operatorname{Mod}}_\Lambda^{#1}}
\def\dim{\operatorname{dim}}

\def\guni{{\Cal G}\operatorname{-Uni}}
\def\gunis{{\Cal G}\operatorname{-Uni}(\SS)}
\def\calp{{\frak p}}
\def\calq{{\frak q}}

\def\BongAdv{{\bf 1}}
\def\GeomIV{{\bf 2}}
\def\GeomI{{\bf 3}}
\def\GeomIII{{\bf 4}}
\def\lebruyn{{\bf 5}}

\topmatter
\title The geometry of uniserial representations of algebras II. Alternate
viewpoints and uniqueness
\endtitle
\rightheadtext{Geometry of uniserial representations. II}

\author Klaus Bongartz and Birge Huisgen-Zimmermann\endauthor

\address FB Mathematik, Universit\"at Gesamthochschule Wuppertal, 42119
Wuppertal, Germany
\endaddress
\email bongartz\@math.uni-wuppertal.de\endemail

\address Department of Mathematics, University of California, Santa Barbara, CA
93106\endaddress
\email birge\@math.ucsb.edu\endemail

 \thanks The research of the second author was partially supported by a
National Science Foundation grant.\endthanks


\abstract  We provide two alternate settings for a family of varieties
modeling the uniserial representations with
fixed sequence of composition factors over a finite dimensional algebra.  The
first is a quasi-projective subvariety of a Grassmannian containing the
members of the mentioned family as a principal affine open cover; among other
benefits, one derives invariance from this intrinsic description.  The second
viewpoint re-interprets the `uniserial varieties' as locally closed
subvarieties of the traditional module varieties; in particular, it exhibits
closedness of the fibres of the canonical maps from the uniserial varieties to
the uniserial representations.      
\endabstract

\endtopmatter
\document

\head 1. Introduction, notation, and further prerequisites \endhead

Let $\Lambda = K\Gamma / I$ be a finite dimensional path algebra modulo
relations, based on a quiver $\Gamma$ and an algebraically closed field $K$ of
coefficients.  In \cite{\GeomI}, a family of affine algebraic varieties was
introduced which not only provides a crisp fit for the interlaced parameters
governing the structure of the uniserial $\la$-modules, but which is concretely
accessible and manageable. In
\cite{\GeomIII, \GeomIV} these `uniserial varieties' were used to derive
representation-theoretic information from quiver and relations.

Originally, the present note was solely aimed at giving a proof for the
invariance of the uniserial varieties under a change of coordinatization of
$\la$  --  invariance up to birational equivalence,
more precisely. It turned out that this goal is most smoothly achieved via a
shift of viewpoint, exhibiting the mentioned family of varieties as an affine
open cover of a quasi-projective subvariety of a Grassmannian.  This
realization somewhat widened the scope of the note.

In fact, we will present two alternate incarnations of the varieties of
\cite{\GeomI} and provide an application for each.  In particular, we will use
a bridge to the classical varieties of representations to give an easy proof
for the fact that the natural maps from the uniserial varieties to the set of
isomorphism types of uniserial $\la$-modules always have closed fibres.  Our
argument is based on the interesting observation that there are no nontrivial
uniserial degenerations of uniserial modules (Section 3). On the side, we
mention a third interpretation of our varieties in terms of Hesselink
stratifications of nullcones, given by L. Le Bruyn \cite{\lebruyn}, which
provides another angle on the subject and potentially lends itself to
generalization. 

For more detail on Section 2, we fix a sequence 
$$\SS = (S(0), \dots, S(l))$$ 
of simple left $\la$-modules $S(i) = \la e(i)/ Je(i)$, where the
$e(i)$ are vertices of $\Gamma$, identified with the corresponding primitive
idempotents of $\la$.  Note that for each uniserial left $\la$-module
$U$ with sequence
$\SS$ of consecutive composition factors, there is at least one path $p$ of
length $l$ passing through the vertices $e(0), e(1), \dots, e(l)$ in this
sequence such that $pU \ne 0$.  Each such path will be referred to as a
{\it mast} of $U$; for short, we will moreover say that $U$ has composition
series $\SS$ and that $p$  is a path passing through
$\SS$.  We briefly recall the definition of the varieties
$V_p$ of
\cite{\GeomI}, where
$p$ traces the paths through
$\SS$.  If $p_0 = e(0), p_1, \dots, p_l$ are the right subpaths of a mast
$p$ of a given uniserial module $U = \la x$, ordered by length (i\.e\., $p =
q_i p_i$, meaning $q_i$ after $p_i$, and
$\len(p_i) = i$), we describe $U$ by recording the effect which multiplication
by arrows has on the basis
$\{p_i x\}$.  Just `a few' of the arising coefficients need to be recorded to
pin down the isomorphism class of $U$:  Given $m \le l$ and an arrow
$\alpha$ in $\Gamma$, we call the pair $(\alpha, p_m)$ a {\it detour} on
$p$  --  and write $(\alpha, p_m) \detour p$  --  in case $\alpha p_m \in
K\Gamma \setminus I$ is different from $p_{m+1}$ but there exists an index
$s > m$ such that $p_s$ ends in the same vertex as $\alpha$.  Letting
$I(\alpha, p_m)$ be the set of all such indices $s$, we then obtain $\alpha
p_m x = \sum_{i \in I(\alpha, p_m)} k_i(\alpha, p_m) p_i x$ for unique scalars
$k_i(\alpha,p_m)$.  It is straightforward to check that the family
$(k_i(\alpha, u))_{i \in I(\alpha,u), (\alpha, u) \detour p}$ determines $U$. 
In
\cite{\GeomI}, it was shown that the set of all such coordinate tuples, as $U$
runs through the uniserials with mast $p$, is an affine algebraic variety, a
defining set of polynomials for which can easily be derived from the quiver
$\Gamma$ and a set of generators for $I$.  The variety $V_p$ comes equipped
with an obvious map $\Phi_p: V_p \rightarrow \{\text{iso types of uniserials
with mast}\ p\}$ given by 
$$k = (k_i(\alpha,u)) \mapsto \la e(0) \bigm/ \biggl( \sum_{(\alpha,u) \detour
p} \la \bigl(\alpha u \ -
\sum_{i \in I(\alpha,u)} k_i(\alpha,u)p_i \bigr) \ \  + \sum_{q \ \text{not a
route on}\ p} \la q \biggr),$$  where a path $q = q e(0)$ is called a {\it
route on} $p$ if the sequence of vertices through which it consecutively passes
is a subsequence of $(e(0), \dots, e(l))$.  One of the virtues of this setup
is the snug fit of the
$V_p$'s; `often', $\Phi_p$ is a bijection, this being for example the case
when $p$ does not pass through an oriented cycle (see \cite{\GeomI}).   

Now $V_{\SS}$ denotes the finite family of irreducible components of the
varieties $V_p$, where $p$ traces the paths through $\SS$.  In case there is
more than one such path, this juxtaposition of partial frames is somewhat
unsatisfactory; it reflects the shortcoming that, up to now, there was no
obvious gluing of these affine patches to a variety incorporating the complete
picture.  This deficiency is removed in Section 2, where we give an
interpretation of the varieties $V_p$ as subvarieties
$\guni(p)$ of a Grassmannian.  While being computationally inaccessible (unless
one takes the detour through the $V_p$'s), the varieties
$\guni(p)$ have the advantage of embedding  --  by dint of their very
construction  --  into a quasiprojective variety
$\gunis$ inside the Grassmannian in which we are operating.  The family of
varieties
$\guni(p)$, where $p$ traces the paths through $\SS$, turns out to be an
affine cover of
$\gunis$ consisting of principal open sets and very natural in the following
sense:  First, there is a canonical map
$\phi_{\SS}$ from $\gunis$ to the set of isomorphism types of uniserials with
composition series $\SS$, the restrictions of which to the patches
$\guni(p)$ are nothing else but the reincarnations of the maps
$\Phi_p$ in the new setting. In other words, if we denote the restriction of
$\phi_{\SS}$ to $\guni(p)$ by $\phi_p$, our isomorphism 
$\guni(p) \rightarrow V_p$ renders the following triangles commutative:

\ignore{
$$\xymatrixcolsep{4pc}\xymatrixrowsep{1pc}
\xy\xymatrix{ V_p \ar[dr]^-{\Phi_p} \ar[dd]_{\cong}\\
 & \save+<23ex,0ex> \drop{\txt{ \{iso types of uniserials in $\lamod$ with mast
$p$\} }} \restore\\
\guni(p) \ar[ur]_-{\phi_p} }\endxy$$ }

\noindent Moreover, the subvarieties $\guni(p)$ of $\gunis$ are unions of
fibres of $\phi_{\SS}$, and so are all intersections $$\guni(p_1) \cap \dots
\cap
\guni(p_r).$$  Thus the variety $\gunis$ provides a joint geometric framework
for the collection of the $V_p$'s which makes the uniqueness statements of
Theorems E and F of \cite{\GeomI} quite obvious in that the definition of
$\gunis$ is intrinsic.  Moreover, this setting will turn out very helpful in
the description of the fibres of the maps $\Phi_p$ in \cite{\GeomIV}.

\head 2. Uniserial varieties as subvarieties of Grassmannians. Uniqueness
\endhead

Let $\SS = (S(0), \dots, S(l))$ be a sequence of simple $\la$-modules with
$S(i) = \la e(i)/ Je(i)$, abbreviate $e(0)$ to $e$, and let $p \in K\Gamma
\setminus I$ be a path through $\SS$.  As before, $p_0,
\dots, p_l$ will be the right subpaths of $p$ with $\len(p_i) = i$, and by
$A_p$ we will denote the $K$-subspace of
$\la e$ spanned by the $p_i$.  Moreover, we set $m = \dim_K \la
e  -  \dim_K A_p = \dim_K \la e - (l+1)$.  Finally, $\Cal G_m(\la e)$ will
stand for the Grassmannian of all $m$-dimensional subspaces of $\la e$.  One
readily checks that the subset $\Cal G_m^{\la} (\la e)$ of $\Cal G_m(\la e)$
consisting of all
$\la$-submodules of $\la e$ having $K$-dimension $m$ is a closed subvariety. 
It is inside this closed subset that the varieties $\guni(p)$ and $\gunis$
will be located.  We start by observing that, given any module $C \in \Cal
G_m^{\la} (\la e)$ with $\la e = A_p \oplus C$, we have $J^i(\la e/ C) = \la
(p_i + C)$, and hence
$\la e/ C$ is uniserial with mast $p$.

\definition{Definition} Let $\guni(p)$ be the set of those objects in
$\Cal G_m^{\la} (\la e)$ which complement the subspace $A_p$ in $\la e$.  We
define $\gunis$ as the union of the sets $\guni(q)$, where $q$ traces the
paths through
$\SS$.

Moreover, we label by $\phi_{\SS}$ the map $$\gunis \longrightarrow
\{\text{iso types of uniserials with composition series}\ \SS \}$$ given by $C
\mapsto \la e/ C$, and by
$\phi_p$ we will denote its restriction to $\guni(p)$.  \enddefinition

By construction, the image of the map $\phi_p$ is just the set of isomorphism
classes of uniserial modules with mast $p$, and hence $\gunis$ equals the set
of all
$\la$-submodules $C$ of $\la e$ with the property that $\la e / C$ is uniserial
with composition series $\SS$.  In particular, the fibres of the restrictions
$\phi_p$ of $\phi_{\SS}$ coincide with fibres of $\phi_{\SS}$, and all
intersections of the patches $\guni(p)$ are unions of fibres.  Moreover, we
have an obvious map
$\psi_p: V_p
\longrightarrow
\guni(p)$ which makes the triangle

\ignore{
$$\xymatrixcolsep{4pc}\xymatrixrowsep{1pc}
\xy\xymatrix{ V_p \ar[dr]^-{\Phi_p} \ar[dd]_{\psi_p} \\
 & \save+<23ex,0ex> \drop{\txt{ \{iso types of uniserials in $\lamod$ with mast
$p$\} }} \restore\\
\guni(p) \ar[ur]_-{\phi_p} }\endxy$$ }

\noindent commutative. Namely, 
$$\psi_p( (k_i(\alpha, u)_{(\alpha, u) \detour p}) = \sum_{(\alpha,u) \detour
p} \la \bigl(\alpha u \ -
\sum_{i \in I(\alpha,u)} k_i(\alpha,u)p_i \bigr) \ \  + \sum_{q \ \text{not a
route on}\ p} \la q $$ satisfies this requirement (see \cite{\GeomI}).

\definition{Example}   Let $\SS = (S(0),S(1))$.  If there are $s \ge 1$ arrows
from
$e(0)$ to $e(1)$, say $\gamma_1, \dots, \gamma_s$, then clearly each
$\la$-submodule $C$ of $\la e(0)$ with the property that $\la e(0) / C$ is
$2$-dimensional with composition series $\SS$ contains the submodule $B = J^2
e(0)\ + \sum_{e'
\ne e(1)} e'Je(0)$, and hence the points of $\gunis$ are in one-to-one
correspondence with the
$(s-1)$-dimensional $K$-subspaces of $J e(0) / B$.  Consequently, $\gunis$ is
a projective variety, a copy of $\Bbb P^{s-1}$ in fact.  Moreover, the masts
$p$ run through the $\gamma_i$, and the
family $(\guni(p))$ provides a canonical affine cover.  

In case $\la = K\Gamma$ is a hereditary algebra, we are presented with a
similarly symmetric situation.  Indeed, if $\SS = (S(0), \dots, S(l))$ is as
before and there are precisely $s(i)$ distinct arrows from the vertex $e(i-1)$
to the vertex $e(i)$, the subset $\gunis$ of $\G_m(\la e(0))$ is a
multi-projective subvariety, namely a copy of
$\Bbb P^{s(1)-1} \times \cdots
\times \Bbb P^{s(l)-1}$.  Again the varieties $V_p$, where $p$ traces the
paths through
$\SS$, constitute the canonical affine cover in the most convenient
coordinatization. 
\enddefinition
  
This illustration suggests the following general connection between the
varieties $V_p$ and the Grassmannians.  The first part is essentially obvious. 
The main point, namely that the $\psi_p$'s are isomorphisms, is equally
elementary; the slightly technical appearance of our argument is due to the
fact that there is no general machinery available for dealing with the
$V_p$'s.

\proclaim{Theorem A}  Each of the sets $\guni(p)$, where $p$ traces the paths
through the sequence $\SS$, is an open affine subvariety of 
$\G_m^{\la}(\la e)$; in fact,
$\guni(p)$ is the intersection of a principal affine subvariety of $\G_m(\la
e)$ with the set of submodules of $\la e$.  Moreover, the maps
$\psi_p$ are isomorphisms $V_p \rightarrow \guni(p)$. \endproclaim

\demo{Proof}  On all counts, it is harmless to assume that there exists a
uniserial module with mast $p$, this being equivalent to $\guni(p)$ being
nonempty and to $V_p$ being nonempty.  In particular, this implies that
$p_0, p_1, \dots, p_l \in \la e$ are $K$-linearly independent in $\la e$
(throughout this proof, we identify paths with their canonical images in
$\la$).  Thus we can supplement the
$p_i$ to a $K$-basis $p_0, p_1, \dots p_l, b_1, \dots, b_m$ of
$\la e$.  Letting $\Cal B$ be an ordered basis of $\bigwedge^m \la e$ recruited
from the subsets of cardinality $m$ of $\{p_i, b_j \}$, with first vector $b_1
\wedge
\dots \wedge b_m$, we then observe that a space $D$ in $\G_m(\la e)$ trivially
intersects $A_p =
\sum_{0 \le i
\le l}K p_i$ if and only if the first of the homogeneous coordinates of
$\bigwedge^m D \in \Bbb P(\bigwedge^m \la e)$ relative to $\Cal B$ is nonzero.

For a proof of the final assertion, start by observing that $\psi_p$ is onto. 
Indeed, each
$\la$-submodule $C$ of $\la e$ which complements the subspace $A_p$ contains
all those paths starting in $e$ which are non-routes on $p$; in fact, $C$ is
generated by the non-routes and some elements of the form $\alpha u - \sum_{i >
\len(u)} k_i p_i$ over $\la$, where
$(\alpha, u)$ is a detour on $p$ and the $k_i$ are scalars; it is thus
sufficient to invoke the definition of $V_p$ to conclude that each object in
$\guni(p)$ is hit by $\psi_p$.  Injectivity of 
$\psi_p$ is obvious in view of the linear independence of the elements $p_i +
\psi_p(k) \in \la e / \psi_p(k)$, $0 \le i \le l$. 

To verify that the map $\psi_p$ and its inverse are morphisms, it will turn out
convenient to choose vectors
$b_1, \dots, b_m$ along the line of the first paragraph, as follows:  First
supplement $p_0, \dots, p_l$ to a basis of
$A_p +
\sum_{(\alpha, u)
\detour p} K\alpha u$ by picking suitable elements $b_1, \dots, b_{m_1}$ from
$\{ \alpha u \mid (\alpha,u) \detour p \}$.  Write $b_i = \alpha_i u_i$ for $i
= 1,\dots, m_1$, and impose the side condition that the sum $\sum_{1 \le i \le
m_1} \len(u_i)$ be maximal.  Next supplement $p_0, \dots, p_l, b_1, \dots,
b_{m_1}$ to a basis of $$A_p \ \ \ + \ 
\sum_{(\alpha, u) \detour p} K\alpha u \ \ \ \  + \sum_{q \in \la e\ \text{a
non-route on}\ p} Kq$$ by picking certain non-routes
$b_{m_1 + 1}, \dots, b_{m_2}$ in $\la e$, and finally supplement
$p_0, \dots, p_l, \allowmathbreak b_1, \dots, b_{m_2}$ to a basis of $\la e$ by
recruiting suitable routes
$b_{m_2 + 1}, \dots, b_{m}$ on $p$.  We claim that each module $C \in \guni(p)$
then has a unique vectorspace basis of the form 
$$\align \{b_i \ \ \  - \sum_{\len(u_i) < j \le l} k_j(\alpha_i, u_i) p_j \mid
1 \le i \le m_1 \} &\cup \{b_{m_1 + 1}, \dots, b_{m_2} \}\\
 &\cup \{b_i \ - \sum_{1 \le j \le l} l_{ij} p_j \mid m_2+1 \le i \le m \},
\endalign$$  where the $k_j(\alpha_i, u_i)$ are scalars and the $l_{ij}$ are
polynomial expressions in the
$k_j(\alpha_i,u_i)$, with polynomials not depending on the choice of $C$.  

Uniqueness is obvious in view of the fact that $C \cap A_p = 0$, whereas
existence follows from the following considerations.  First note that, apart
from the claim concerning the $l_{ij}$, existence is once more an immediate
consequence of the fact that $\la e = A_p
\oplus C$; just keep in mind that all non-routes on $p$ are zero modulo $C$ and
that $\alpha u + C$ belongs to $J^{\len(u) + 1}(\la e / C)$ for any detour
$(\alpha,u)$, whereas each residue class
$p_i + C$ lies outside $J^{i+1}(\la e / C)$.  So our task is reduced to showing
that, given any route $v$ on $p$, there exist polynomials $\calq_r$ such that,
for any choice of
$C$ in $\guni(p)$, the element $v \in \la e$ is congruent to $\sum_{1 \le r\le
l} \calq_r (k_j(\alpha_i, u_i)) p_r$ modulo $C$; here the scalars
$k_j(\alpha_i, u_i)$ depend of course on $C$.  We will do this by induction on
$l - d$, where $d$ is maximal with the property that $p_d$ is a right subpath
of $v$. In case $l - d = 0$, we have $v = p$.  So suppose this difference is
positive.  It is clearly harmless to assume that
$v$ is not a right subpath of $p$, i\.e\.,
$v = v'\beta p_d$, where
$(\beta, p_d)$ is a detour on $p$ and $\len(v') \ge 0$.  By our
length-maximal choice of $b_1,
\dots, b_{m_1}$, there are scalars $h_r, k_s$ such that 
$\beta p_d = \sum_{1 \le r \le l} h_r p_r + \sum_{1 \le s \le m_1} k_s \alpha_s
u_s$  with $\len(u_s) \ge d$ whenever $k_s \ne 0$.  Again, since $p_r + C
\notin J^{r+1}(\la e/ C)$, we see that, for all $r$ with
$h_r \ne 0$, we have $\len(p_r) \ge d+1$.   To compute $v + C = v' \beta p_d
+ C$, let us multiply the equation for $\beta p_d$ by $v'$ from the left and
individually focus on the residue classes $v' p_r + C$ and $v'\alpha_s u_s + C$
with $r \ge d+1$ and $\len(u_s) \ge d$, respectively.  If the path
$v' p_r$ is not a route on
$p$, its residue class modulo $C$ is always zero.  In the remaining case,
our induction hyposthesis yields polynomials
$\calp_t$, $1 \le t \le l$, depending only on $v'$ and $r$ such that $v' p_r +
C =  \sum_{1 \le t \le l}
\calp_t(k_j(\alpha_i,u_i)) p_t + C$.  As for the other terms:  By the
definition of the scalars
$k_z(\alpha_s,u_s)$ we have $v'\alpha_s u_s + C = \sum_{\len(u_s) < z \le l}
k_z(\alpha_s,u_s) v' p_z + C$, and in view of $\len(u_s) \ge d$, we obtain $v'
p_z + C =  \sum_{1 \le t \le l}
\calp_t(k_j(\alpha_i,u_i)) v'p_t + C$ for the terms of the latter sum.  This
completes the proof of our intermediate claim.      
   
Note that the argument of the preceding paragraph holds for arbitrary routes
$v$ on $p$ so, in particular, for all routes $\alpha u$, where $(\alpha,u)$ is
a detour.  It follows that
$V_p$ is isomorphic to its projection onto the components
$(k_j(\alpha_i,u_i))_{j \in I(\alpha_i, u_i), 1 \le i \le m_1}$.  Again we use
our special basis for $\la e$ to form an ordered basis
$\Cal B$ for $\bigwedge^m \la e$ consisting of
$m$-element subsets and containing $b_1 \wedge \dots \wedge b_m$ as first
vector.  The first paragraph of the proof permits us to identify the points of
$\guni(p)$ with points of $\bigwedge^m \la e$ having first coordinate $1$,
which identifies 
$\psi_p$ with the map sending a point $k = (k_i(\alpha,u))$ to $c_1\wedge \dots
\wedge c_m$, where $c_1, \dots, c_m$ is the basis of the space $\psi_p(k)$
specified in the preceding claim; in particular, we have $c_i = b_i -
\sum_{\len(u_i) < j \le l} k_j(\alpha_i,u_i) p_j$ for $1 \le i \le m_1$.  The
claim clearly guarantees that
$\psi_p$ is a morphism.  As for
$\psi_p^{-1}$, it will be enough to show that, up to signs, all the scalars
$k_j(\alpha_i,u_i)$ show up as coefficients in the expansion of $c_1 \wedge
\dots \wedge c_m$ relative to $\Cal B$.  But since, up to a factor $\pm 1$,
the elements $d_{ij} = d_{ij1}\wedge
\dots \wedge d_{ijm}$, $1 \le i \le m_1$, with $d_{ijr} = b_r$ for $r \ne i$ and
$d_{iji} = p_j$ belong to $\Cal B$, this is clear. \qed               
\enddemo

\proclaim{Corollary B} $\gunis$ is a quasiprojective variety, and the family
$(\guni(p))$, where $p$ traces the paths through $\SS$, is an affine cover
consisting of principal open subsets of $\gunis$. \qed\endproclaim 

The next consequence of our observation proves Theorems E and F of
\cite{\GeomI}. 

\proclaim{Corollary C} If $\la \cong \la'$ and $\SS'$ is the sequence of
simple left 
$\la'$-modules corresponding to $\SS$ under an algebra isomorphism, the
varieties $\gunis$ and $\guni(\SS')$ are isomorphic.  

In particular, given a coordinatization $K \Gamma'/I'$ of $\la'$, there is a
one-to-one correspondence between the collection of all irreducible components
of the
$V_{p'}$, where
$p'$ traces the paths in $\Gamma'$ through $\SS'$, and the irreducible
components of the $V_p$, where p traces the paths in $\Gamma$ through
$\SS$, such that the partners under this pairing are birationally equivalent.
\endproclaim

\demo{Proof} Let $f: \la \rightarrow \la'$ be an algebra isomorphism and $e' =
f(e)$, then clearly the induced isomorphism $G_m(\la e) \cong G_m(\la' e')$ 
restricts to an isomorphism $\gunis \cong \guni(\SS')$.  Since the affine
patches $\guni(p')$ of
$\guni(\SS')$ relative to any coordinatization of $\la'$ form an open cover of
$\guni(\SS')$, the second claim is a consequence of the first. \qed \enddemo

\head 3. The uniserial varieties as subvarieties of the classical ones. 
Closedness of the fibres of the canonical maps. \endhead

Let $e_1, \dots, e_n$ be the distinct vertices of $\Gamma$, and let $\Gamma*$
be the union of $\{e_1, \dots, e_n\}$ with the set of arrows of $\Gamma$. 
Recall that, given
$d \in \Bbb N$, the variety $\modla{d}$ of (bounden) representations of $\la
= K\Gamma / I$ with $K$-dimension $d$ consists of those points $$x =
(\alpha(x)) \in \prod_{\alpha \in \Gamma*}
\Hom_K(K^{d_{s(\alpha)}}, K^{d_{t(\alpha)}}),$$ the components of which satisfy
the relations in
$I$.  By $R$ we denote the canonical map from
$\modla{d}$ to the set of isomorphism classes of $d$-dimensional left
$\la$-modules.  If $\SS = (S(0), \dots, S(l))$ is a sequence of simple modules
as before and $d = l+1$, the subset $\rep \SS \subseteq
\modla{d}$ of the points corresponding to uniserial modules with
composition series $\SS$ is an open subvariety of
$\modla{d}$ which is closed under the standard action of $\GL_d$;  indeed, if
$p$ traces the paths through $\SS$, then $\rep
\SS$  is the union of the open sets consisting of the points in $\modla{d}$
corresponding to the uniserials with mast $p$ (i\.e\., the points
$x$ with $\gamma_l(x) \cdots \gamma_1(x) \ne 0$, provided that $p = \gamma_l
\cdots \gamma_1$).  It is within these open sets that we will find closed
subvarieties
$\rep(p)$ which will again turn out to be canonically isomorphic to the $V_p$.

In the sequel, we will write the points of $x \in \modla{d}$ as sequences of
$d \times d$-matrices, representing the maps $\alpha(x)$ relative to the
canonical basis of $K^d$. 

\definition{Definition} Given a path $p = \gamma_l \cdots \gamma_1$ through
$\SS$, where $\gamma_i$ is an arrow
$e(i-1) \rightarrow e(i)$, we define
$\rep(p)$ to be the closed subset of $\rep \SS$ consisting of those points $x$
which are represented by tuples of strictly lower triangular matrices
$(\alpha(x))_{\alpha \in \Gamma*}$ having the additional property that the
$i$-th column of the matrix $\gamma_i(x)$ equals the $(i+1)$-st canonical basis
vector for $1 \le i \le l$. \enddefinition

It is clear that the image $R(\rep(p))$ is the full set of isomorphism classes
of uniserial modules with mast $p$; in other words, the image of $R$
restricted to $\rep(p)$ coincides with the image of $\Phi_p$.

\proclaim{Theorem D} There is a canonical isomorphism $\rho_p: \rep(p)
\rightarrow V_p$ of varieties which makes the following diagram commutative:

\ignore{
$$\xymatrixcolsep{4pc}\xymatrixrowsep{1pc}
\xy\xymatrix{ V_p \ar[dr]^{\Phi_p}\\
 & \save+<23ex,0ex> \drop{\txt{ \{iso types of uniserials in $\lamod$ with mast
$p$\} }} \restore\\
\rep(p) \ar[ur]_{R} \ar[uu]^{\rho_p}_{\cong} }\endxy$$ }

\endproclaim

\demo{Proof}  Define $\rho_p (x)$ to be $(k_i^x (\alpha, p_m))_{i \in
I(\alpha,p_m), (\alpha,p_m)
\detour p}$, where $k_i^x(\alpha, p_m)$ is the the entry in position
$(i,m)$ of the matrix $\alpha(x)$.  In other words, $\rho_p$ amounts to a
projection onto a certain subset of coordinates of $\modla{d}$ and is thus a
morphism.  All other entries of $x$ are either $0$ or $1$, depending solely on
$p$, which shows that $\rho_p$ is `essentially' equal to the identity: Indeed,
a point $x \in \rep(p)$ corresponding to a uniserial module
$U$ essentially just consists of matrix representations of the multiplications
by arrows, relative to a suitable basis of the form $(p_0 z, \dots, p_l
z)$, $z$ a top element of $U$, coinciding with the canonical basis of
$K^n$.  Focusing on a matrix $\alpha(x)$ for some arrow
$\alpha$, we note that, for a
non-detour $(\alpha,p_m)$, the entry in position $(i,m)$  is
$1$ or $0$ depending on whether or not $i=m+1$ and $\alpha p_m = p_{m+1}$; in
case of a detour, on the other hand, the corresponding entry is zero for all $i
\notin I(\alpha,p_m)$.  This clearly guarantees that $\rho_p$ is an
isomorphism. 

Commutativity of the diagram follows from the obvious fact that $x$ and
$\rho_p(x)$ describe the same uniserial module with mast $p$, up to
isomorphism. \qed
\enddemo 
\medskip

\noindent{\it Remark.}  Again let $d = l+1$, where $l = \len(p)$.  As will
become clear in \cite{\GeomIV}, the variety $\rep (p)$ is  --  via fibre
bundles and geometric quotients  --  intrinsically related to the open
subvariety of $\modla{d}$ consisting of all those representations which are
not annihilated by $p$ (note that these are automatically uniserial). 
Before we had grown aware of the interpretation in terms of Grassmannians
described in Section 2, we made use of this relationship to prove
Corollary C.  Our present viewpoint has a twofold edge over the previous
one, however:  First, it yields an obvious action of a unipotent group on
$\rep(p)$ such that the orbits coincide with the isomorphism classes of
uniserials (see \cite{\GeomIV}); second, it provides us with a global
quasi-projective variety which glues the various $V_p$'s together to a
variety representing {\it all} uniserial modules with a fixed sequence of
composition factors.   
\smallskip

As an application, we give an elementary proof for the fact that the fibres of
the canonical maps $\Phi_p$ are closed.  In view of Theorems A and D, this is
clearly tantamount to closedness of the fibres of the sibling maps $\phi_p$, as
well as to closedness inside $\rep \SS$ of the intersections of the
$\GL_d$-orbits of $\modla{d}$ with $\rep \SS$.  

\proclaim{Proposition E} Suppose that $U,U'\in \lamod$ are non-isomorphic
uniserial modules of the same $K$-dimension. Then $U$ does not degenerate
to $U'$.\endproclaim

\demo{Proof} The beginning of our induction on $\dim_K U$ is trivial since, in
the variety of representations, the points corresponding to semisimple modules
have closed orbits. Now suppose that $\dim_K U > 1$, and assume that $U'$ is a
degeneration of $U$. As is well known, we then have $\dim_K
\Hom_\la(U,X) \le \dim_K \Hom_\la(U',X)$ and $\dim_K
\Hom_\la(X,U) \le \dim_K \Hom_\la(X,U')$ for all $X\in
\lamod$ (see, e.g., \cite{\BongAdv}). If $S$ denotes the socle of $U$, then
plugging $S=X$ into the second inequality, we see that $S$ is the socle of $U'$
as well. Theorem 2.4 of \cite{\BongAdv} now tells us that
$Q= U/S$ degenerates to $Q'= U'/S$, and hence the induction hypothesis entails
$Q\cong Q'$. Since $U$ and $U'$ are non-isomorphic but have the same
$K$-dimension, we therefore obtain an isomorphism
$\Hom_\la(U',U)\cong \Hom_\la(Q,U)$, which yields $\dim_K
\Hom_\la(U,U)> \dim_K \Hom_\la(Q,U)= \dim_K \Hom_\la(U',U)$. This contradicts
the first of the cited inequalities. \qed\enddemo

\proclaim{Corollary F} The fibres of $\Phi_p: V_p \rightarrow \{\text{iso types
of uniserials with mast}\ p\}$ are closed subvarieties of $V_p$.
$\qed$ \endproclaim

\enddocument